\input amstex
\documentstyle{amsppt}
\font\tenscr=rsfs10 \font\sevenscr=rsfs7 \font\fivescr=rsfs5
\skewchar\tenscr='177 \skewchar\sevenscr='177
\skewchar\fivescr='177
\newfam\scrfam \textfont\scrfam=\tenscr \scriptfont\scrfam=\sevenscr
\scriptscriptfont\scrfam=\fivescr
\define\scr#1{{\fam\scrfam#1}}
%%%%%%%%%%%%%%%%%%%%%%%%%%%%%%%%%%%%%%%%%%%%%%%%%
%Version: Alex Degtyarev (home) 28/09/98  16:11:30
\def\stydate{May 10, 2002}

\chardef\tempcat\catcode`\@
%\ifx\undefined\amstexloaded
%\input amstex \else\catcode`\@\tempcat\fi \expandafter\ifx\csname
%amsppt.sty\endcsname\relax\input amsppt.sty \fi
\let\tempcat\undefined

\immediate\write16{This is LABEL.DEF by A.Degtyarev <\stydate>}
\expandafter\ifx\csname label.def\endcsname\relax\else
  \message{[already loaded]}\endinput\fi
\expandafter\edef\csname label.def\endcsname{%
  \catcode`\noexpand\@\the\catcode`\@\edef\noexpand\styname{LABEL.DEF}%
  \def\expandafter\noexpand\csname label.def\endcsname{\stydate}%
    \toks0{}\toks2{}}
\catcode`\@11
\def\labelmesg@ {LABEL.DEF: }
{\edef\temp{\the\everyjob\W@{\labelmesg@<\stydate>}}
\global\everyjob\expandafter{\temp}}

%%%%%%%%%%%%%%%%%%%% Handy stuff
\def\@car#1#2\@nil{#1}
\def\@cdr#1#2\@nil{#2}
\def\eat@bs{\expandafter\eat@\string}
\def\eat@ii#1#2{}
\def\eat@iii#1#2#3{}
\def\eat@iv#1#2#3#4{}
\def\@DO#1#2\@{\expandafter#1\csname\eat@bs#2\endcsname}
\def\@N#1\@{\csname\eat@bs#1\endcsname}
\def\@Nx{\@DO\noexpand}
\def\@Name#1\@{\if\@undefined#1\@\else\@N#1\@\fi}
\def\@Ndef{\@DO\def}
\def\@Ngdef{\global\@Ndef}
\def\@Nedef{\@DO\edef}
\def\@Nxdef{\global\@Nedef}
\def\@Nlet{\@DO\let}
\def\@undefined#1\@{\@DO\ifx#1\@\relax\true@\else\false@\fi}
\def\@@addto#1#2{{\toks@\expandafter{#1#2}\xdef#1{\the\toks@}}}
\def\@@addparm#1#2{{\toks@\expandafter{#1{##1}#2}%
    \edef#1{\gdef\noexpand#1####1{\the\toks@}}#1}}
\def\make@letter{\edef\t@mpcat{\catcode`\@\the\catcode`\@}\catcode`\@11 }
\def\donext@{\expandafter\egroup\next@}
\def\x@notempty#1{\expandafter\notempty\expandafter{#1}}
%%%%% Defines #1 to be #2 in lower case
\def\lc@def#1#2{\edef#1{#2}%
    \lowercase\expandafter{\expandafter\edef\expandafter#1\expandafter{#1}}}
%%%%% Finds a text in a comma separated list
\newif\iffound@
\def\find@#1\in#2{\found@false
    \DNii@{\ifx\next\@nil\let\next\eat@\else\let\next\nextiv@\fi\next}%
    \edef\nextiii@{#1}\def\nextiv@##1,{%
    \edef\next{##1}\ifx\nextiii@\next\found@true\fi\FN@\nextii@}%
    \expandafter\nextiv@#2,\@nil}
%%%%% Disable \outer'ness
{\let\head\relax\let\specialhead\relax\let\subhead\relax
\let\subsubhead\relax\let\proclaim\relax
\gdef\let@relax{\let\head\relax\let\specialhead\relax\let\subhead\relax
    \let\subsubhead\relax\let\proclaim\relax}}
%%%%% Hacks
\newskip\@savsk
% We add a tiny skip in order to make \@esphack work correct
\let\@ignorespaces\ignorespaces
\def\@ignorespacesp{\ifhmode
  \ifdim\lastskip>\z@\else\penalty\@M\hskip-1sp%
        \penalty\@M\hskip1sp \fi\fi\@ignorespaces}
\def\ignorespaces{\protect\@ignorespacesp}
% However, we have to redefine some control sequences
%%??\def~{\unskip\nobreak\ \@ignorespaces}
\def\@bsphack{\relax\ifmmode\else\@savsk\lastskip
  \ifhmode\edef\@sf{\spacefactor\the\spacefactor}\fi\fi}
\def\@esphack{\relax
  \ifx\penalty@\penalty\else\penalty\@M\fi   % if this is after \nobreak
  \ifmmode\else\ifhmode\@sf{}\ifdim\@savsk>\z@\@ignorespacesp\fi\fi\fi}
%%%%% Check for a *, \nofrills, and []
\let\@frills@\identity@
\let\@txtopt@\identyty@
\newif\if@star
\newif\if@write\@writetrue
\def\@numopt@{\if@star\expandafter\eat@\fi}
\def\checkstar@#1{\DN@{\@writetrue
  \ifx\next*\DN@####1{\@startrue\checkstar@@{#1}}%
      \else\DN@{\@starfalse#1}\fi\next@}\FN@\next@}
\def\checkstar@@#1{\DN@{%
  \ifx\next*\DN@####1{\@writefalse#1}%
      \else\DN@{\@writetrue#1}\fi\next@}\FN@\next@}
\def\checkfrills@#1{\DN@{%
  \ifx\next\nofrills\DN@####1{#1}\def\@frills@####1{####1\nofrills}%
      \else\DN@{#1}\let\@frills@\identity@\fi\next@}\FN@\next@}
\def\checkbrack@#1{\DN@{%
    \ifx\next[\DN@[####1]{\def\@txtopt@########1{####1}#1}%
    \else\DN@{\let\@txtopt@\identity@#1}\fi\next@}\FN@\next@}
\def\check@therstyle#1#2{\bgroup\DN@{#1}\ifx\@txtopt@\identity@\else
        \DNii@##1\@therstyle{}\def\@therstyle{\DN@{#2}\nextii@}%
    \expandafter\expandafter\expandafter\nextii@\@txtopt@\@therstyle.\@therstyle
    \fi\donext@}

%%%%%%%%%%%%%%%%%%%% \@input, \include, .aux files, etc.
\newread\@inputcheck
\def\@input#1{\openin\@inputcheck #1 \ifeof\@inputcheck \W@
  {No file `#1'.}\else\closein\@inputcheck \relax\input #1 \fi}

\def\loadstyle#1{\edef\next{#1}%
    \DN@##1.##2\@nil{\if\notempty{##2}\else\def\next{##1.sty}\fi}%
    \expandafter\next@\next.\@nil\lc@def\next@\next
    \expandafter\ifx\csname\next@\endcsname\relax\input\next\fi}

\let\pagebody@\pagebody
\let\pagetop@\empty
\let\pagebot@\empty
\let\@Xend\empty
\def\pagebody{\pagetop@\pagebody@\pagebot@\@Xend}
\let\@Xclose\empty

\newwrite\@Xmain
\newwrite\@Xsub
\def\W@X{\write\@Xout}
\def\make@Xmain{\global\let\@Xout\@Xmain\global\let\end\endmain@
  \xdef\@Xname{\jobname}\xdef\@inputname{\jobname}}
\begingroup
\catcode`\(\the\catcode`\{\catcode`\{12
\catcode`\)\the\catcode`\}\catcode`\}12
\gdef\W@count#1((\lc@def\@tempa(#1)%
    \def\\##1(\W@X(\global##1\the##1))%
    \edef\@tempa(\W@X(%
        \string\expandafter\gdef\string\csname\space\@tempa\string\endcsname{)%
        \\\pageno\\\cnt@toc\\\cnt@idx\\\cnt@glo\\\footmarkcount@
        \@Xclose\W@X(}))\expandafter)\@tempa)
\endgroup
\def\readaux{\bgroup\checkbrack@\readaux@}
\let\begin\readaux
\def\readaux@{%
    \W@{>>> \labelmesg@ Run this file twice to get x-references right}%
    \global\everypar{}%
    {\def\\{\global\let}%
        \def\/##1##2{\gdef##1{\wrn@command##1##2}}%
        \disablepreambule@cs}%
    \make@Xmain{\make@letter\setboxz@h{\@input{\@txtopt@{\@Xname.aux}}%
            \lc@def\@tempa\jobname\@Name\open@\@tempa\@}}%
  \immediate\openout\@Xout\@Xname.aux%
    \immediate\W@X{\relax}\egroup}
\everypar{\global\everypar{}\readaux}
{\toks@\expandafter{\topmatter}
\global\edef\topmatter{\noexpand\readaux\the\toks@}}
\let\@@end@@\end

\def\@Xclose@{{\def\@Xend{\ifnum\insertpenalties=\z@
        \W@count{close@\@Xname}\closeout\@Xout\fi}%
    \vfill\supereject}}
\def\endmain@{\@Xclose@
    \W@{>>> \labelmesg@ Run this file twice to get x-references right}%
    \@@end@@}
\def\disablepreambule@cs{\\\disablepreambule@cs\relax}

\def\include#1{\bgroup
  \ifx\@Xout\@Xsub\DN@{\errmessage
        {\labelmesg@ Only one level of \string\include\space is supported}}%
    \else\edef\@tempb{#1}\clearpage
      \DN@##1 {\if\notempty{##1}\edef\@tempb{##1}\DN@####1\eat@ {}\fi\next@}%
    \DNii@##1.{\edef\@tempa{##1}\DN@####1\eat@.{}\next@}%
        \expandafter\next@\@tempb\eat@{} \eat@{} %
    \expandafter\nextii@\@tempb.\eat@.%
        \relaxnext@
      \if\x@notempty\@tempa
          \edef\nextii@{\write\@Xmain{%
            \noexpand\string\noexpand\@input{\@tempa.aux}}}\nextii@
        \ifx\undefined\@includelist\found@true\else
                    \find@\@tempa\in\@includelist\fi
            \iffound@\ifx\undefined\@noincllist\found@false\else
                    \find@\@tempb\in\@noincllist\fi\else\found@true\fi
            \iffound@\lc@def\@tempa\@tempa
                \if\@undefined\close@\@tempa\@\else\edef\next@{\@Nx\close@\@tempa\@}\fi
            \else\xdef\@Xname{\@tempa}\xdef\@inputname{\@tempb}%
                \W@count{open@\@Xname}\global\let\@Xout\@Xsub
            \openout\@Xout\@tempa.aux \W@X{\relax}%
            \DN@{\let\end\endinput\@input\@inputname
                    \@Xclose@\make@Xmain}\fi\fi\fi
  \donext@}
\def\includeonly#1{\edef\@includelist{#1}}
\def\noinclude#1{\edef\@noincllist{#1}}

%%%%%%%%%%%%%%%%%%%% Numeric styles
\def\arabicnum#1{\number#1}

\def\Romannum#1{\expandafter\uppercase\expandafter{\romannumeral#1}}
\def\alphnum#1{\ifcase#1\or a\or b\or c\or d\else\@ialph{#1}\fi}
\def\@ialph#1{\ifcase#1\or \or \or \or \or e\or f\or g\or h\or i\or j\or
    k\or l\or m\or n\or o\or p\or q\or r\or s\or t\or u\or v\or w\or x\or y\or
    z\else\fi}
\def\Alphnum#1{\ifcase#1\or A\or B\or C\or D\else\@Ialph{#1}\fi}
\def\@Ialph#1{\ifcase#1\or \or \or \or \or E\or F\or G\or H\or I\or J\or
    K\or L\or M\or N\or O\or P\or Q\or R\or S\or T\or U\or V\or W\or X\or Y\or
    Z\else\fi}

%%%%%%%%%%%%%%%%%%%% Counters
\def\ST@P{step}
\def\ST@LE{style}
\def\N@M{no}
\def\F@NT{font@}
%   #1 is the name of the counter to be defined,
%   #2 is the counter this one depends upon,
\outer\def\newcounter{\checkbrack@{\expandafter\newcounter@\@txtopt@{{}}}}
{\let\newcount\relax
\gdef\newcounter@#1#2#3{{%
    \toks@@\expandafter{\csname\eat@bs#2\N@M\endcsname}%
    \DN@{\alloc@0\count\countdef\insc@unt}%
    \ifx\@txtopt@\identity@\expandafter\next@\the\toks@@
        \else\if\notempty{#1}\global\@Nlet#2\N@M\@#1\fi\fi
    \@Nxdef\the\eat@bs#2\@{\if\@undefined\the\eat@bs#3\@\else
            \@Nx\the\eat@bs#3\@.\fi\noexpand\arabicnum\the\toks@@}%
  \@Nxdef#2\ST@P\@{}%
  \if\@undefined#3\ST@P\@\else
    \edef\next@{\noexpand\@@addto\@Nx#3\ST@P\@{%
             \global\@Nx#2\N@M\@\z@\@Nx#2\ST@P\@}}\next@\fi
    \expandafter\@@addto\expandafter\@Xclose\expandafter
        {\expandafter\\\the\toks@@}}}}
\outer\def\copycounter#1#2{%
    \@Nxdef#1\N@M\@{\@Nx#2\N@M\@}%
    \@Nxdef#1\ST@P\@{\@Nx#2\ST@P\@}%
    \@Nxdef\the\eat@bs#1\@{\@Nx\the\eat@bs#2\@}}
\outer\def\everystep{\checkstar@\everystep@}
\def\everystep@#1{\if@star\let\next@\gdef\else\let\next@\@@addto\fi
    \@DO\next@#1\ST@P\@}
%   #1 is the counter whose style is to be changed,
\def\counterstyle#1{\@Ngdef\the\eat@bs#1\@}
\def\advancecounter#1#2{\@N#1\ST@P\@\global\advance\@N#1\N@M\@#2}
\def\setcounter#1#2{\@N#1\ST@P\@\global\@N#1\N@M\@#2}
\def\counter#1{\refstepcounter#1\printcounter#1}
\def\printcounter#1{\@N\the\eat@bs#1\@}
\def\refcounter#1{\xdef\@lastmark{\printcounter#1}}
\def\stepcounter#1{\advancecounter#1\@ne}
\def\refstepcounter#1{\stepcounter#1\refcounter#1}
\def\savecounter#1{\@Nedef#1@sav\@{\global\@N#1\N@M\@\the\@N#1\N@M\@}}
\def\restorecounter#1{\@Name#1@sav\@}

%%%%%%%%%%%%%%%%%%%% Warnings
\def\warning#1#2{\W@{Warning: #1 on input line #2}}
\def\warning@#1{\warning{#1}{\the\inputlineno}}
\def\wrn@@Protect#1#2{\warning@{\string\Protect\string#1\space ignored}}
\def\wrn@@label#1#2{\warning{label `#1' multiply defined}{#2}}
\def\wrn@@ref#1#2{\warning@{label `#1' undefined}}
\def\wrn@@cite#1#2{\warning@{citation `#1' undefined}}
\def\wrn@@command#1#2{\warning@{Preamble command \string#1\space ignored}#2}
\def\wrn@@option#1#2{\warning@{Option \string#1\string#2\space is not supported}}
\def\wrn@@reference#1#2{\W@{Reference `#1' on input line \the\inputlineno}}
\def\wrn@@citation#1#2{\W@{Citation `#1' on input line \the\inputlineno}}
\let\wrn@reference\eat@ii
\let\wrn@citation\eat@ii
%% Disable wornings: works with
%% \Protect, \label, \ref, \cite, \command, \option, \reference, \citation
%% (with FONT.DEF)   \font
%% (with DEBUG.DEF)  \nocite
\def\nowarning#1{\if\@undefined\wrn@\eat@bs#1\@\wrn@option\nowarning#1\else
        \@Nlet\wrn@\eat@bs#1\@\eat@ii\fi}
\def\printwarning#1{\if\@undefined\wrn@@\eat@bs#1\@\wrn@option\printwarning#1\else
        \@Nlet\wrn@\eat@bs#1\expandafter\@\csname wrn@@\eat@bs#1\endcsname\fi}
\printwarning\Protect \printwarning\label \printwarning\ref
\printwarning\cite \printwarning\command \printwarning\option

%%%%%%%%%%%%%%%%%%%% Hyperrefs
{\catcode`\#=12\gdef\@lH{#}}
\def\@@HREF#1{}
\def\@HREF#1#2{\@@HREF{a #1}{\let\@@HREF\eat@#2}\@@HREF{/a}}
\def\@@Hf#1{file:#1} \let\@Hf\@@Hf
\def\@@Hl#1{\if\notempty{#1}\@lH#1\fi} \let\@Hl\@@Hl
\def\@@Hname#1{\@HREF{name="#1"}{}} \let\@Hname\@@Hname
\def\@@Href#1{\@HREF{href="#1"}} \let\@Href\@@Href
\ifx\undefined\pdfoutput
  \csname newcount\endcsname\pdfoutput
\else
  \def\pdflinkattr{attr{/C [0 0.9 0.9]}}
  \let\pdflinkbegin\empty
  \let\pdflinkend\empty
  \def\@pdfHf#1{file {#1}}
  \def\@pdfHl#1{name {#1}}
  \def\@pdfHname#1{\pdfdest name{#1}xyz\relax}
  \def\@pdfHref#1#2{\pdfstartlink \pdflinkattr goto #1\relax
    \pdflinkbegin#2\pdflinkend\pdfendlink}
  \def\@ifpdf#1#2{\ifnum\pdfoutput>\z@\expandafter#1\else\expandafter#2\fi}
  \def\@Hf{\@ifpdf\@pdfHf\@@Hf}
  \def\@Hl{\@ifpdf\@pdfHl\@@Hl}
  \def\@Hname{\@ifpdf\@pdfHname\@@Hname}
  \def\@Href{\@ifpdf\@pdfHref\@@Href}
\fi
\def\@Hr#1#2{\if\notempty{#1}\@Hf{#1}\fi\@Hl{#2}}
\def\@localHref#1{\@Href{\@Hr{}{#1}}}
\def\@countlast#1{\@N#1last\@}
\def\@@countref#1#2{\global\advance#2\@ne
  \@Nxdef#2last\@{\the#2}\@tocHname{#1\@countlast#2}}
\def\@countref#1{\@DO\@@countref#1@HR\@#1}

%\Href@-xxx#1#2#3
\def\Href@@#1{\@N\Href@-#1\@}
\def\Href@#1#2{\@N\Href@-#1\@{\@Hl{@#1-#2}}}
%\Hname@-xxx#1
\def\Hname@#1{\@N\Hname@-#1\@}
%\Hlast@-xxx
\def\Hlast@#1{\@N\Hlast@-#1\@}
\def\cntref@#1{\global\@DO\advance\cnt@#1\@\@ne
  \@Nxdef\Hlast@-#1\@{\@DO\the\cnt@#1\@}\Hname@{#1}{@#1-\Hlast@{#1}}}
\def\HyperRefs#1{\global\@Nlet\Hlast@-#1\@\empty
  \global\@Nlet\Hname@-#1\@\@Hname
  \global\@Nlet\Href@-#1\@\@Href}
\def\NoHyperRefs#1{\global\@Nlet\Hlast@-#1\@\empty
  \global\@Nlet\Hname@-#1\@\eat@
  \global\@Nlet\Href@-#1\@\eat@}

%% Hyperrefs for label/ref
\HyperRefs{label} {\catcode`\-11
\gdef\@labelref#1{\Hname@-label{r@-#1}}
\gdef\@xHref#1{\Href@-label{\@Hl{r@-#1}}} }
%% Hyperrefs for toc
\HyperRefs{toc}
\def\@HR#1{\if\notempty{#1}\string\@HR{\Hlast@{toc}}{#1}\else{}\fi}

%\def\HyperRefs{\gdef\@labelref##1{\@Hname{r@-##1}}%
% \gdef\@xHref##1{\@localHref{r@-##1}}}
%\def\NoHyperRefs{\global\let\@labelref\eat@\global\let\@xHref\eat@}
%\def\TOCRefs{\global\let\@tocHname\@Hname\global\let\@tocHref\@localHref}
%\def\NoTOCRefs{\global\let\@tocHname\eat@\global\let\@tocHref\eat@}
%%\def\@HR#1{\if\notempty{#1}\string\@HR{\cnt@toclast}{#1}\else{}\fi}
%\HyperRefs
%%\TOCRefs

%%\def\@CLR#1{\special{color push #1}}
%%\def\@eCLR{\special{color pop}}

%%%%%%%%%%%%%%%%%%%% Labels, x-references
\def\bftext{\ifmmode\fam\bffam\else\bf\fi}
\let\@lastmark\empty
\let\@lastlabel\empty
\def\lastmark{\@lastmark}
\let\lastlabel\empty
\let\everylabel\relax
\let\everylabel@\eat@
\let\everyref\relax
\def\newlabel{\bgroup\everylabel\newlabel@}
\def\newlabel@#1#2#3{\if\@undefined\r@-#1\@\else\wrn@label{#1}{#3}\fi
  {\let\protect\noexpand\@Nxdef\r@-#1\@{#2}}\egroup}
\def\w@ref{\bgroup\everyref\w@@ref}
\def\w@@ref#1#2#3#4{%
  \if\@undefined\r@-#1\@{\bftext??}#2{#1}{}\else%
   \@xHref{#1}{\@DO{\expandafter\expandafter#3}\r@-#1\@\@nil}\fi
  #4{#1}{}\egroup}%\null}}
\def\@@@xref#1{\w@ref{#1}\wrn@ref\@car\wrn@reference}
\def\@xref#1{\rom{\@@@xref{#1}}}
\let\xref\@xref
\def\pageref#1{\w@ref{#1}\wrn@ref\@cdr\wrn@reference}
\def\thepage{\ifnum\pageno<\z@\romannumeral-\pageno\else\number\pageno\fi}
\def\label@{\@bsphack\bgroup\everylabel\label@@}
\def\label@@#1#2{\everylabel@{{#1}{#2}}%
  \@labelref{#2}%
  \let\thepage\relax
  \def\protect{\noexpand\noexpand\noexpand}%
  \edef\@tempa{\edef\noexpand\@lastlabel{#1}%
    \W@X{\string\newlabel{#2}{{\@lastmark}{\thepage}}{\the\inputlineno}}}%
  \expandafter\egroup\@tempa\@esphack}
\def\label#1{\label@{#1}{#1}}
\def\fn@P@{\relaxnext@
    \DN@{\ifx[\next\DN@[####1]{}\else
        \ifx"\next\DN@"####1"{}\else\DN@{}\fi\fi\next@}%
    \FN@\next@}
\def\eat@fn#1{\ifx#1[\expandafter\eat@br\else
  \ifx#1"\expandafter\expandafter\expandafter\eat@qu\fi\fi}
\def\eat@br#1]#2{}
\def\eat@qu#1"#2{}
{\catcode`\~\active\lccode`\~`\@
\lowercase{\global\let\@@P@~\gdef~{\protect\@@P@}}}
\def\Protect@@#1{\def#1{\protect#1}}
\def\disable@special{\let\W@X@\eat@iii\let\label\eat@
    \def\footnotemark{\protect\fn@P@}%
  \let\footnotetext\eat@fn\let\footnote\eat@fn
    \let\refcounter\eat@\let\savecounter\eat@\let\restorecounter\eat@
    \let\advancecounter\eat@ii\let\setcounter\eat@ii
  \let\ifvmode\iffalse\Protect@@\@@@xref\Protect@@\pageref\Protect@@\nofrills
    \Protect@@\\\Protect@@~}
\let\notoctext\identity@
\def\W@X@#1#2#3{\@bsphack{\disable@special\let\notoctext\eat@
    \def\chapter{\protect\chapter@toc}\let\thepage\relax
    \def\protect{\noexpand\noexpand\noexpand}#1%
  \edef\next@{\if\@undefined#2\@\else\write#2{#3}\fi}\expandafter}\next@
    \@esphack}
\newcount\cnt@toc
\def\writeauxline#1#2#3{\W@X@{\cntref@{toc}\let\tocref\@HR}
  \@Xout{\string\@Xline{#1}{#2}{#3}{\thepage}}}
{\let\newwrite\relax
\gdef\@openin#1{\make@letter\@input{\jobname.#1}\t@mpcat}
\gdef\@openout#1{\global\expandafter\newwrite\csname
tf@-#1\endcsname
   \immediate\openout\@N\tf@-#1\@\jobname.#1\relax}}
\def\@@openout#1{\@openout{#1}%
  \@@addto\readaux@{\immediate\closeout\@N\tf@-#1\@}}
\def\auxlinedef#1{\@Ndef\do@-#1\@}
\def\@Xline#1{\if\@undefined\do@-#1\@\expandafter\eat@iii\else
    \@DO\expandafter\do@-#1\@\fi}
\def\beginW@{\bgroup\def\do##1{\catcode`##112 }\dospecials\do\@\do\"
    \catcode`\{\@ne\catcode`\}\tw@\immediate\write\@N}
\def\endW@toc#1#2#3{{\string\tocline{#1}{#2\string\page{#3}}}\egroup}
\def\do@tocline#1{%
%%  The file version
    \if\@undefined\tf@-#1\@\expandafter\eat@iii\else
        \beginW@\tf@-#1\@\expandafter\endW@toc\fi
%%  The \toks version
%       \if\@undefined\the#1@@\@\else
%           \global\addto{\the#1@}{\tocline{#2}{#3\page{#4}}}\fi
} \auxlinedef{toc}{\do@tocline{toc}}

\let\protect\empty
\def\Protect#1{\if\@undefined#1@P@\@\PROTECT#1\else\wrn@Protect#1\empty\fi}
\def\PROTECT#1{\@Nlet#1@P@\@#1\edef#1{\noexpand\protect\@Nx#1@P@\@}}
\def\pdef#1{\edef#1{\noexpand\protect\@Nx#1@P@\@}\@Ndef#1@P@\@}

\Protect\operatorname \Protect\operatornamewithlimits
\Protect\qopname@ \Protect\qopnamewl@ \Protect\text
\Protect\topsmash \Protect\botsmash \Protect\smash
\Protect\widetilde \Protect\widehat \Protect\thetag
\Protect\therosteritem
% Fonts:
\Protect\Cal \Protect\Bbb \Protect\bold \Protect\slanted
\Protect\roman \Protect\italic \Protect\boldkey
\Protect\boldsymbol \Protect\frak \Protect\goth \Protect\dots
% Symbols
\Protect\cong \Protect\lbrace \let\{\lbrace \Protect\rbrace
\let\}\rbrace
\let\root@P@@\root \def\root@P@#1{\root@P@@#1\of}
\def\root#1\of{\protect\root@P@{#1}}

\def\frills{\ignorespaces\@txtopt@}
\def\frillsnotempty#1{\x@notempty{\@txtopt@{#1}}}
\def\numberline{\@numopt@}
\newif\if@theorem
\let\@therstyle\eat@
\def\@headtext@#1#2{{\disable@special\let\protect\noexpand
    \def\chapter{\protect\chapter@rh}%
    \edef\next@{\noexpand\@frills@\noexpand#1{#2}}\expandafter}\next@}
\let\AmSrighthead@\rightheadtext
\def\rightheadtext{\checkfrills@{\@headtext@\AmSrighthead@}}
\let\AmSlefthead@\leftheadtext
\def\leftheadtext{\checkfrills@{\@headtext@\AmSlefthead@}}
% #1 refers to the style,
% #2 refers to the style in the toc,
% #3 refers to the counter,
% #4 represents the AmSTeX's counterpart (we use \end... because
%   of \outer-ness), and
% #5 is the text to be typeset
\def\@head@@#1#2#3#4#5{\@Name\pre\eat@bs#1\@\if@theorem\else
    \@frills@{\csname\expandafter\eat@iv\string#4\endcsname}\relax
        \ifx\protect\empty\@N#1\F@NT\@\fi\fi
    \@N#1\ST@LE\@{\counter#3}{#5}%
  \if@write\writeauxline{toc}{\eat@bs#1}{#2{\counter#3}\@HR{#5}}\fi
    \if@theorem\else\expandafter#4\fi
    \ifx#4\endhead\ifx\@txtopt@\identity@\else
        \headmark{\@N#1\ST@LE\@{\counter#3}{\frills\empty}}\fi\fi
    \@Name\post\eat@bs#1\@\ignorespaces}
\ifx\undefined\endhead\Invalid@\endhead\fi
\def\@head@#1{\checkstar@{\checkfrills@{\checkbrack@{\@head@@#1}}}}
% #1 is the name,
% #2 is the counter, and
% #3 is the title text
\def\@thm@@#1#2#3{\@Name\pre\eat@bs#1\@
    \@frills@{\csname\expandafter\eat@iv\string#3\endcsname}
    {\@theoremtrue\check@therstyle{\@N#1\ST@LE\@}\frills
            {\counter#2}\@theoremfalse}%
    \@DO\envir@stack\end\eat@bs#1\@
    \@N#1\F@NT\@\@Name\post\eat@bs#1\@\ignorespaces}
\def\@thm@#1{\checkstar@{\checkfrills@{\checkbrack@{\@thm@@#1}}}}
% #1 is the name,
% #2 is the counter,
% #3 is the name of the corresponding table (lof, lot, etc.)
% #4 is either \topcaption or \botcaption
%   #5 is the caption text
\def\@capt@@#1#2#3#4#5\endcaption{\bgroup
    \edef\@tempb{\global\footmarkcount@\the\footmarkcount@
    \global\@N#2\N@M\@\the\@N#2\N@M\@}%
    \def\shortcaption##1{\global\def\sh@rtt@xt####1{##1}}\let\sh@rtt@xt\identity@
%    \let\notoctext\identity@
%%% To handle \nofrills !!!
%   \DN@##1##2##3{\false@\fi\iftrue}%
%   \ifx\@frills@\identity@\else\let\notempty\next@\fi
    \DN@{#4{\@tempb\@N#1\ST@LE\@{\counter#2}}}%
    \if\notempty{#5}\DNii@{\next@\@N#1\F@NT\@}\else\let\nextii@\next@\fi
    \nextii@#5\endcaption
  \if@write\writeauxline{#3}{\eat@bs#1}{{} \@HR{\@N#1\ST@LE\@{\counter#2}%
    \if\notempty{#5}.\enspace\fi\sh@rtt@xt{#5}}}\fi
  \global\let\sh@rtt@xt\undefined\egroup}
\def\@capt@#1{\checkstar@{\checkfrills@{\checkbrack@{\@capt@@#1}}}}
\let\captiontextfont@\empty

\ifx\undefined\subsubheadfont@\def\subsubheadfont@{\it}\fi
\ifx\undefined\proclaimfont\def\proclaimfont{\sl}\fi
\ifx\undefined\proclaimfont@\let\proclaimfont@\proclaimfont\fi
\def\proclaimfont{\proclaimfont@}
\ifx\undefined\definitionfont@\def\AmSdeffont@{\rm}
    \else\let\AmSdeffont@\definitionfont@\fi
\ifx\undefined\remarkfont@\def\remarkfont@{\rm}\fi

\def\newfont@def#1#2{\if\@undefined#1\F@NT\@
    \@Nxdef#1\F@NT\@{\@Nx.\expandafter\eat@iv\string#2\F@NT\@}\fi}
% #1 is the name (and the style),
% #2 is the style in the toc,
% #3 is the counter, and
% #4 is the AmSTeX's counterpart to be used (in the \end... form):
\def\newhead@#1#2#3#4{{%
    \gdef#1{\@therstyle\@therstyle\@head@{#1#2#3#4}}\newfont@def#1#4%
    \if\@undefined#1\ST@LE\@\@Ngdef#1\ST@LE\@{\headstyle}\fi
    \if\@undefined#2\@\gdef#2{\headtocstyle}\fi
  \@@addto\moretocdefs@{\\#1#1#4}}}
\outer\def\newhead#1{\checkbrack@{\expandafter\newhead@\expandafter
    #1\@txtopt@\headtocstyle}}
% #1 is the default title (like Theorem, Lemma, etc.),
% #2 is the name to be defined,
% #3 refers to the counter (which should be defined separately),
% #4 is the AmSTeX's counterpart: \endproclaim, \endremark, \endAmSdef
\outer\def\newtheorem#1#2#3#4{{%
    \gdef#2{\@thm@{#2#3#4}}\newfont@def#2#4%
    \@Nxdef\end\eat@bs#2\@{\noexpand\revert@envir
        \@Nx\end\eat@bs#2\@\noexpand#4}%
  \if\@undefined#2\ST@LE\@\@Ngdef#2\ST@LE\@{\proclaimstyle{#1}}\fi}}%
% #1 is the default title (like Figure, Table, etc.),
% #2 is the name to be defined,
% #3 refers to the counter (which should be defined separately),
% #4 is the name of the corresponding table (toc, lof, lot, etc.)
% #5 is either \topcaption or \botcaption
\outer\def\newcaption#1#2#3#4#5{{\let#2\relax
  \edef\@tempa{\gdef#2####1\@Nx\end\eat@bs#2\@}%
    \@tempa{\@capt@{#2#3{#4}#5}##1\endcaption}\newfont@def#2\endcaptiontext%
  \if\@undefined#2\ST@LE\@\@Ngdef#2\ST@LE\@{\captionstyle{#1}}\fi
  \@@addto\moretocdefs@{\\#2#2\endcaption}\newtoc{#4}}}
{
% #1 is the name of the table (toc, lof, lot, etc.)
\outer\gdef\newtoc#1{{%
    \@DO\ifx\do@-#1\@\relax
%%  The \toks version
%    \global\expandafter\newtoks\csname the#4@@\endcsname
    \global\auxlinedef{#1}{\do@tocline{#1}}{}%
    \@@addto\tocsections@{\make@toc{#1}{}}\fi}}}

\toks@\expandafter{\itembox@}
\toks@@{\bgroup\let\therosteritem\identity@\let\rm\empty
  \let\@Href\eat@\let\@Hname\eat@
  \edef\next@{\edef\noexpand\@lastmark{\therosteritem@}}\donext@}
\edef\itembox@{\the\toks@@\the\toks@}
\def\firstitem@false{\let\iffirstitem@\iffalse
    \global\let\lastlabel\@lastlabel}

\let\rosteritemrefform\therosteritem
\let\rosteritemrefseparator\empty
\def\rosteritemref#1{\hbox{\rosteritemrefform{\@@@xref{#1}}}}
\def\local#1{\label@\@lastlabel{\lastlabel-i#1}}

\def\xRef@P@{\gdef\lastlabel}
\def\xRef#1{\@xref{#1}\protect\xRef@P@{#1}}

\def\iref@P@{\gdef\lastref}
\def\itemref#1#2{\rosteritemref{#1-i#2}\protect\iref@P@{#1}}
\def\iref#1{\@xref{#1}\rosteritemrefseparator\itemref{#1}}

\def\eqtag{\tag\counter\equation}
\def\eqref#1{\thetag{\@@@xref{#1}}}
\def\tagform@#1{\ifmmode\hbox{\rm\else\rom{\fi
        (\ignorespaces#1\unskip)\iftrue}\else}\fi}

\let\AmSfnote@\makefootnote@
\def\makefootnote@#1{\bgroup\let\footmarkform@\identity@
  \edef\next@{\edef\noexpand\@lastmark{#1}}\donext@\AmSfnote@{#1}}

\def\clearpage{\ifnum\insertpenalties>0\line{}\fi\vfill\supereject}

\def\proof{\checkfrills@{\checkbrack@{%
    \check@therstyle{\@frills@{\demo}{\frills{Proof}}{}}
        {\frills{}\envir@stack\endremark\envir@stack\enddemo}%
  \envir@stack\endproof\ignorespaces}}}
\def\everyendproof{\qed}
\def\endproof{\nofrillscheck{\frills@{\everyendproof}\revert@envir\endproof\enddemo}}

\let\AmSref\ref
\let\AmSrefstyle\refstyle
\let\plaincite\cite
\def\citei@#1,{\citeii@#1\eat@,}
\def\citeii@#1\eat@{\w@ref{#1}\wrn@cite\@car\wrn@citation}
\def\mcite@#1;{\plaincite{\citei@#1\eat@,\unskip}\mcite@i}
\def\mcite@i#1;{\DN@{#1}\ifx\next@\endmcite@
  \else, \plaincite{\citei@#1\eat@,\unskip}\expandafter\mcite@i\fi}
\def\endmcite@{\endmcite@}
\def\cite#1{\mcite@#1;\endmcite@;}
\PROTECT\cite
\def\refstyle#1{\AmSrefstyle{#1}\uppercase{%
    \ifx#1A\relax \def\@ref@##1{\AmSref\xdef\@lastmark{##1}\key##1}%
    \else\ifx#1C\relax \def\@ref@##1{\AmSref\no\counter\refno}%
        \else\def\@ref@{\AmSref}\fi\fi}}
\refstyle A
\newcounter\refno\null
\newif\ifRefs
\gdef\Refs{\checkstar@{\checkbrack@{\csname AmSRefs\endcsname
  \nofrills{\frills{References}%
  \if@write\writeauxline{toc}{vartocline}{\@HR{\frills{References}}}\fi}%
  \def\ref{\@ref@}\Refstrue\ignorespaces}}}
\let\ref\xref

\newif\iftoc
\pdef\tocbreak{\iftoc\hfil\break\fi}
\def\tocsections@{\make@toc{toc}{}}
\let\moretocdefs@\empty
\def\newtocline@#1#2#3{%
  \edef#1{\def\@Nx#2line\@####1{\@Nx.\expandafter\eat@iv
        \string#3\@####1\noexpand#3}}%
  \@Nedef\no\eat@bs#1\@{\let\@Nx#2line\@\noexpand\eat@}%
    \@N\no\eat@bs#1\@}
\def\MakeToc#1{\@@openout{#1}}
\def\newtocline#1#2#3{\Err@{\Invalid@@\string\newtocline}}
\def\make@toc#1#2{\penaltyandskip@{-200}\aboveheadskip
    \if\notempty{#2}
        \centerline{\headfont@\ignorespaces#2\unskip}\nobreak
    \vskip\belowheadskip \fi
%%  The file version
    \@openin{#1}\relax
%%  The \toks version
%   \if\@undefined\the#1@@\@
%       \the\@N\the#1@@\@\global\@N\the#1@@\@\fi
    \vskip\z@}
\def\contents{\readaux\checkfrills@{\checkbrack@{\@contents@}}}
\def\@contents@{\toc@{\frills{Contents}}\envir@stack\endcontents%
    \def\nopagenumbers{\let\page\eat@}\let\newtocline\newtocline@\toctrue
  \def\@HR{\Href@{toc}}%
  \def\tocline##1{\csname##1line\endcsname}
  \edef\caption##1\endcaption{\expandafter\noexpand
    \csname head\endcsname##1\noexpand\endhead}%
    \ifmonograph@\def\vartoclineline{\Chapterline}%
        \else\def\vartoclineline##1{\sectionline{{} ##1}}\fi
  \let\\\newtocline@\moretocdefs@
    \ifx\@frills@\identity@\def\\##1##2##3{##1}\moretocdefs@
        \else\let\tocsections@\relax\fi
    \def\\{\unskip\space\ignorespaces}\let\maketoc\make@toc}
\def\endcontents{\tocsections@\vskip-\lastskip\revert@envir\endcontents
    \endtoc}

% \selectf@nt is for future extensions (like RUSSIAN.TEX)
\if\@undefined\selectf@nt\@\let\selectf@nt\identity@\fi
\def\Err@math#1{\Err@{Use \string#1\space only in text}}
\def\textonlyfont@#1#2{%
    \def#1{\RIfM@\Err@math#1\else\edef\f@ntsh@pe{\string#1}\selectf@nt#2\fi}%
    \PROTECT#1}
\tenpoint

% #1 is the name of the switch to be defined
%   #2 is the default font switch
\def\newshapeswitch#1#2{\gdef#1{\selectsh@pe#1#2}\PROTECT#1}
% #1 is the name of the switch
% #2 is the current shape
% #3 is the shape to be used
\def\shapeswitch#1#2#3{\@Ngdef#1\string#2\@{#3}}
% These shapes are used by \rom
\shapeswitch\rm\bf\bf \shapeswitch\rm\tt\tt
\shapeswitch\rm\smc\smc
\newshapeswitch\em\it
% These shapes are used by \em and \emph
\shapeswitch\em\it\rm \shapeswitch\em\sl\rm
\def\selectsh@pe#1#2{\relax\if\@undefined#1\f@ntsh@pe\@#2\else
    \@N#1\f@ntsh@pe\@\fi}

\def\@itcorr@{\leavevmode
    \edef\prevskip@{\ifdim\lastskip=\z@ \else\hskip\the\lastskip\relax\fi}\unskip
    \edef\prevpenalty@{\ifnum\lastpenalty=\z@ \else
        \penalty\the\lastpenalty\relax\fi}\unpenalty
    \/\prevpenalty@\prevskip@}
\def\rom@P@#1{\@itcorr@{\selectsh@pe\rm\rm#1}}
\def\rom{\protect\rom@P@}
%%%%    \Rom will unconditionally switch to \rm, like \rom in AmSppt
\def\Rom@P@#1{\@itcorr@{\rm#1}}
\def\Rom{\protect\Rom@P@}
{\catcode`\-11 \HyperRefs{idx} \HyperRefs{glo}
\newcount\cnt@idx \global\cnt@idx=10000
\newcount\cnt@glo \global\cnt@glo=10000
\gdef\writeindex#1{\W@X@{\cntref@{idx}}\tf@-idx
 {\string\indexentry{#1}{\Hlast@{idx}}{\thepage}}}
\gdef\writeglossary#1{\W@X@{\cntref@{glo}}\tf@-glo
 {\string\glossaryentry{#1}{\Hlast@{glo}}{\thepage}}}
}
\def\emph#1{\@itcorr@\bgroup\em\ignorespaces#1\unskip\egroup
  \DN@{\DN@{}\ifx\next.\else\ifx\next,\else\DN@{\/}\fi\fi\next@}\FN@\next@}
\def\makequoteactive{\catcode`\"\active}
{\makequoteactive\gdef"{\FN@\quote@}
\gdef\quote@{\ifx"\next\DN@"##1""{\quoteii{##1}}\else\DN@##1"{\quotei{##1}}\fi\next@}}
\let\quotei\eat@
\let\quoteii\eat@
\def\MakeIndex{\@openout{idx}}
\def\MakeGlossary{\@openout{glo}}

%%%%%%%%%%%%%%%%%%%%%%  Just helpful things %%%%%%%%%%%%%%%%%%%%%%%%%%%%
\def\endofpar#1{\ifmmode\ifinner\endofpar@{#1}\else\eqno{#1}\fi
    \else\leavevmode\endofpar@{#1}\fi}
\def\endofpar@#1{\unskip\penalty\z@\null\hfil\hbox{#1}\hfilneg\penalty\@M}

\newdimen\normalparindent\normalparindent\parindent
\def\firstparindent#1{\everypar\expandafter{\the\everypar
  \global\parindent\normalparindent\global\everypar{}}\parindent#1\relax}

%% Commands to disable
\@@addto\disablepreambule@cs{%
    \\\readaux\relax
    \\\begin\relax
    \\\readaux@\relax
    \\\@openout\eat@
    \\\@@openout\eat@
    \/\Monograph\empty
    \/\MakeIndex\empty
    \/\MakeGlossary\empty
    \/\MakeToc\eat@
    \/\HyperRefs\eat@
    \/\NoHyperRefs\eat@
}

\csname label.def\endcsname

%%%%%%%%%%%%%%%%%%%%%%%% Definitions for my papers %%%%%%%%%%%%%%%%%%%%%%%

\def\punct#1#2{\if\notempty{#2}#1\fi}
\def\sppunct{\punct{.\enspace}}
\def\varpunct#1#2{\if\frillsnotempty{#2}#1\fi}

\def\headstyle#1#2{\numberline{#1\sppunct{#2}}\ignorespaces#2\unskip}
\def\headtocstyle#1#2{\numberline{#1\punct.{#2}}\space #2}

\def\specialtocstyle#1#2{#2}
\newcounter\section\null
\newcounter\subsection\section
\newcounter\subsubsection\subsection
\newhead\specialsection[\specialtocstyle]\null\endspecialhead
\newhead\section\section\endhead
\newhead\subsection\subsection\endsubhead
\newhead\subsubsection\subsubsection\endsubsubhead
\def\firstappendix{\global\sectionno0 %
  \counterstyle\section{\Alphnum\sectionno}%
    \global\let\firstappendix\empty}

\def\appendixtocstyle#1#2{\space\numberline{Appendix #1\sppunct{#2}}#2}
\newhead\appendix[\appendixtocstyle]\section\endhead

\let\endAmSdef\enddefinition
\def\proclaimstyle#1#2{\numberline{#2\varpunct{.\enspace}{#1}}\frills{#1}}
\copycounter\thm\subsubsection
%\newcounter[\subsubsectionno]\thm\subsection
\theorem\thm\endproclaim
\proposition\thm\endproclaim
\lemma\thm\endproclaim
\corollary\thm\endproclaim
\definition\thm\endAmSdef
\example\thm\endAmSdef

\def\captionstyle#1#2{\frills{#1}\numberline{\varpunct{ }{#1}#2}}
\newcounter\figure\null
\newcounter\table\null
\newcaption{Figure}\figure\figure{lof}\botcaption
\newcaption{Table}\table\table{lot}\topcaption

\copycounter\equation\subsubsection
%\newcounter[\subsubsectionno]\equation\subsection

%\endinput

%%%%%%%%%%%%%%%%%%%%%%%%%%%%%%%%%%%%%%%%%%%%%%%%%%

\def\C{{\Bbb C}}
\def\R{{\Bbb R}}
\def\Z{{\Bbb Z}}

\def\NN{\scr N}
\def\kk{\varkappa}
\def\GG#1{\widetilde G_2(\R^{#1+2})}
\def\EE#1{\widetilde E_2(\R^{#1+2})}
\def\gg#1{\widetilde\tau_{2,#1+2}}
\def\nGG#1{ G_2(\R^{#1+2})}

\def\ngg#1{\tau_{2,#1+2}}

\nologo\NoBlackBoxes

 \topmatter
\title
Abundance of real lines on real projective hypersurfaces
\endtitle
\author S.~Finashin, V.~Kharlamov
\endauthor
\address Middle East Technical University,
Department of Mathematics\endgraf Ankara 06531 Turkey
\endaddress
%\email  serge@metu.edu.tr \endemail
\address
Universit\'{e} de Strasbourg et IRMA (CNRS)\endgraf 7 rue Ren\'{e}
Descartes 67084 Strasbourg Cedex, France
\endaddress
%\affil \endaffil
%\address \endorsees
%\email \endemail
%\date \enddate
\thanks The second author was partially funded by the
ANR-09-BLAN-0039-01 grant of {\it Agence Nationale de la
Recherche}, and is a member of FRG Collaborative Research "Mirror
Symmetry \& Tropical Geometry" (Award No. $0854989$).
\endthanks
%\dedicatory \enddedicatory
%\keywords Enumerative geometry, real lines, projective hypersurfaces\endkeywords
\subjclassyear{2000} \subjclass 14N15,
14P25, 14J80 \endsubjclass
%\abstract
%We show that a generic real projective n-dimensional hypersurface
%of degree (2n-1) contains "many" real lines, namely, not less than
%(2n-1)!!, which is approximately  the square root of the number of
%complex lines. This estimate is based on the interpretation of a
%suitable signed count of the lines as the Euler number of an
%appropriate bundle.
%\endabstract

\endtopmatter

\hskip1.8in There is geometry in the humming of the strings.

\hskip1.8in There is music in the spacing of the spheres.

\hskip2.5in  PYTHAGORAS (6th century BC)
 \vskip6mm

\document

\section{Revelations}
Our aim is to show that in the case of a generic real hypersurface
$X$ of degree $2n-1$ in a projective space of dimension $n+1$
 the number $\NN_\R$ of real lines on $X$
is not less than approximately  the square root of the number
$\NN_\C$ of complex lines.  More precisely, $\NN_\R\ge (2n-1)!!,$
while due to Don Zagier \cite{GM} $\NN_\C\sim
\sqrt{\frac{27}{\pi}}(2n-1)^{2n-\frac32}$, so that
$$
\log \NN_\R\ge n\log 2n+ O(n)=\frac12\log \NN_\C.
$$

Note that $\NN_\R$, unlike $\NN_\C$, depends not only on $n$ but
also on the choice of $X$. The key point of our estimate is an
appropriate signed count of the real lines that makes the sum
invariant. This sum, which we denote by $\NN^e_\R$, is nothing but
the Euler number of a suitable vector bundle (see \ref{Euler}).
Its evaluation gives finally the conclusion:
$$
  \NN_\R\ge \NN^e_\R=(2n-1)!!.\eqtag\label{bound}
$$

We came to it in 2010--2011, when we started working on a survey
on real cubic hypersurfaces. Realizing that Segre's division of
real lines on cubic surfaces in two species, elliptic and
hyperbolic, leads to a remarkable relation $h-e=3$  between the
number of hyperbolic lines, $h$, and number $e$ of elliptic, we
looked for conceptual explanations and generalizations.

Later, we
recollected that the non-vanishing of the Stiefel-Whitney number
of the same bundle (which implies existence of real lines) was
proved much earlier
  by O.~Debarre and L.~Manivel in \cite{DM}, and that the
Euler number was used by J.~Solomon in his thesis \cite{Sol} to
perform a signed count of real lines on quintic 3-folds. In
September 2011 we were informed by C.~Okonek and A.~Teleman that
they obtained the bound (\ref{bound}) and that they followed a
similar approach. It is when an account of their proof was posted
on arXiv that we decided to compose this note and to present our
proof (in our opinion, a bit  simpler and more geometric) and a
few related observations, like disclosing precise relation to
Welschinger's indices, see \ref{Welschinger}, and revealing the
``secret'' of the ``magic formula'' $h-e=3$, see \ref{magic}.

\section{Origin of Symmetry}

\subsection{The polar  correspondence}
The Grassmannian $\GG{n}$ of oriented 2-planes can be canonically
identified with the purely imaginary quadric $Q^n(\C)$ defined in
$P^{n+1}(\C)$ by equation $x_0^2+x_1^2+\dots+x_{n+1}^2=0$. Namely,
a real oriented plane $p\subset\R^{n+2}$ defines a real oriented
line $\ell_{\R}\subset P^{n+1}(\R)$, which decomposes its
complexification $\ell_{\C}\subset P^{n+1}(\C)$ into a pair of
halves permuted by the complex conjugation. Each half intersects
$Q^n(\C)$ at a single point. We give privilege to that half,
$h_l$, whose complex orientation bounds the given orientation of
$\ell_{\R}$. The map $p\mapsto q=Q^n(\C)\cap h_l$ yields a
diffeomorphism
 $\kk\:\GG{n}\to Q^n(\C)$,
 which we call
{\it  the polar correspondence}.

In the reverse direction, the polar correspondence can be
described as follows. A point $q\in Q^n(\C)\subset P^{n+1}(\C)$ is
represented by a complex line $L_q\subset\C^{n+2}$, the complex
line $L_q$ projects to $\R^{n+2}$ into a real plane $p$ by
averaging $v\mapsto \frac12 (v+\bar v)$, and the complex
orientation of $L_q$ is pushed forward to orient the plane $p$.

 This yields the following interpretation for the tautological
 bundles of the oriented Grassmannian, $\GG{n}$, and of the
 non-orientable one, $G_2(\R^{n+2})$.

\proposition\label{taut}The polar correspondence $\GG{n}\to
Q^n(\C)$ identifies the tautological (oriented 2-plane)
bundle
$\gg{n}\:\EE{n}\to\GG{n}$ with the restriction, $\tau_{n} :
E^n(\C)\to Q^n(\C)$, to $Q^n(\C)$ of the tautological (complex
line) bundle on $ P^{n+1}(\C)$. This yields identification of the
tautological bundle $\ngg{n}$ over $G_2(\R^{n+2})$ with the
quotient of $\tau_n$ by the complex conjugation.
\qed\endproposition

\subsection{Characteristic classes}
The polar correspondence leads to a simple explicit formula for
the Euler class, $e(Sym^{2m-1}(F))$, of symmetric powers of the dual
tautological bundles, $F=\gg{n}^*$. Namely, the polar
correspondence identifies $F$ with the real 2-bundle underlying
complex line bundle $L=\varkappa^*(\tau_n^*)$, so,
$F\otimes\C=L\oplus \bar L$, and $Sym^k(F)\otimes\C=Sym^k(
F\otimes\C)=L^k\oplus L^{k-1}\bar L\dots \oplus \bar L^k$, for any
$k\ge 1$. The complex conjugation interchanges $L^a\bar L^b$ with
$L^b\bar L^a$, thus,
 for any odd $k=2m-1$,
 the real vector bundle $Sym^{2m-1}(F)$
is  the real part of  $L^{2m-1}\oplus L^{2m-2}\bar
L\oplus\dots\oplus \bar L^{2m-1}$. We define an isomorphism
between $L^{2m-1}\oplus L^{2m-2}\bar L\oplus\dots \oplus L^{m}\bar
L^{m-1}$ and $Sym^{2m-1}( F)$ by projecting
 $v\mapsto \frac12 (v+\bar v)$.
 It transports the complex orientation of
 $L^{2m-1}\oplus L^{2m-2}\bar L\oplus\dots \oplus L^{m}\bar L^{m-1}$ to an orientation of
$Sym^{2m-1}(F)$, so that the Euler class of the latter becomes
well defined and equal to the Chern class of the former one.

\proposition\label{classe} Under the above orientation convention,
$$e(Sym^{2m-1}(F))= c_{m}(L^{2m-1}\oplus L^{2m-2}\bar L\oplus\dots \oplus
L^{m}\bar L^{m-1})= (2m-1)!!c_1(L)^m. \qed
$$
\endproposition

\section{Absolution}

\subsection{Euler's numbers}\label{Euler}
  By Proposition \ref{taut}, the polar correspondence identifies the
total space
  $E$ of  the bundle $Sym^{2n-1}(\ngg{n}^*)$
with the quotient by complex conjugation of the total space of a
complex $n$-bundle $(\tau^*_n)^{2n-1}\oplus( \tau_n^*)^{2n-2}(\bar
\tau^*_n)\oplus\dots \oplus( \tau_n^*)^{n}(\bar \tau_n^*)^{n-1}$
over a complex n-manifold, $Q^n(\C)$. Since any anti-holomorphic
involution on a complex manifold of even complex dimension
preserves the orientation, the manifold $E$ inherits an
orientation.

Let $\NN^e_\R$ denote the Euler number of $Sym^{2n-1}(\ngg{n}^*)$
with respect to this orientation.
 Given a section $s\:\nGG{n}\to E$ and an isolated zero $x$ of
 $s$, let $I^e_x(s)$ denote its {\it index}, that is the local Euler number of $s$ at $x$.
 If $s$ has only isolated zeros, then
$\NN^e_\R=\sum_x I^e_x(s)$, and we let $\NN_\R(s)=\sum_x\vert
 I^e_x(s)\vert$.

\theorem\label{intermediate} \roster\item $\NN^e_\R=(2n-1)!!$.
\item
For any section $s$ having only isolated zeros,
$\NN_\R(s)\ge (2n-1)!! $.
\endroster
\endtheorem

\proof By definition, twice $\NN_\R^e$ equals  the Euler number of
$Sym^{2n-1}(\gg{n}^*)$,  which in accordance with Proposition
\ref{classe} is equal to $(2n-1)!!c_1(\tau_n^*)^n[Q^n(\C)]$, while
$c_1(\tau_n^*)^n[Q^n(\C)]=\deg Q^n=2$. This gives (1), wherefrom
(2) is straightforward.
\endproof

Consider a real homogeneous polynomial $f\in H^0(P^{n+1}, \Cal
O(d))$ defining a real hypersurface $X$ of degree $d=2n-1$ in
$P^{n+1}$.
  Restricting $f$ to planes $p\subset\R^{n+2}$,
we obtain a section $s_f$ of $Sym^{2n-1}(\ngg{n}^*)$. Its zeros
represent real lines lying on $X$. We call a real line
$\ell\subset X$ {\it isolated} if it represents an isolated zero
of $s_f$, and attribute to it {\it the real multiplicity} equal,
by definition, to $\vert I^e_\ell(s_f)\vert$. Thus, Theorem
\ref{intermediate} implies the following result.

 \corollary\label{estimate} Assume that all the real lines $\ell\subset X$
of a real hypersurface $X\subset P^{n+1}$ of degree $2n-1$ are
isolated. Then $\sum_{\ell\subset X}\vert I^e_\ell(s_f)\vert\ge
(2n-1)!!$. \qed\endcorollary

As is known (see, for example, \cite{BVV}), for a generic  $f$
(that is for a generic $X$)  the section $s_f$ is transverse to
the zero section, and so, all the lines $\ell\subset X$ are
isolated and all have indices $I_\ell^e(s_f)=\pm1$.
 Let $\NN^+_\R$ and $\NN^-_\R$ denote the number of real lines of indices $1$ and $-1$, respectively.
 Thus, for a generic $f$, the number of real lines is
 $$\NN^+_\R+\NN^-_\R\ge\NN^+_\R-\NN^-_\R=\NN_\R^e=(2n-1)!!$$

\corollary For a generic hypersurface $X\subset P^{n+1}$ of degree
$(2n-1)$, the number of real lines $\ell\subset X$ is finite and
bounded from below by $(2n-1)!!$. \qed\endcorollary

\remark{Remark} Our choice of an orientation of $E$ and of the
multiplicities of $\ell\subset X$ is a bit different from that in
\cite{OT}. The orientation is not essential for us, since $|\NN_\R^e|$
and $|I^e_x(s)|$ remain invariant, while the multiplicities make a
difference. In fact, the complex multiplicities used in \cite{OT}
are bounded from below by
 our indices $|I^e_x(s)|$.
\endremark

\section{Showbiz}\label{showbiz}
\subsection{Cubic surfaces: 3 versus 27}
Recall, that any nonsingular cubic surface contains exactly
$27=c_4( Sym^3(\tau^*_{2,4}(\C)))[G_2(\C^4)]$ complex lines (here
$\tau_{2,4}(\C)$ denotes the tautological bundle over the complex
Grassmannian $G_2(\C^4)$), and it implies, in particular, that for
a nonsingular cubic surface, a section $s_f$, which can be defined
for $Sym^3(\tau^*_{2,4}(\C))$ as well as for  $
Sym^3(\tau^*_{2,4})$, is in the both cases transversal to the zero
section. As a result, Theorem \ref{intermediate} reads as identity
$\NN^e_\R=\NN^+_\R- \NN^-_\R=3$ and estimate $\NN_\R\ge3$.

The estimate $\NN_\R\ge 3$ was known already to Schl\"afli
\cite{Sch}, who proved that the number of real lines on a
non-singular real cubic surface takes only the values 27, 15, 7,
and 3. Segre \cite{Segre} was probably the first to observe that
the real lines on a nonsingular real cubic surface can be divided
in two species, named by him hyperbolic and elliptic, and that
numerically the distribution between the species, $\NN_\R=h+e$, is
as follows: 27=15+12, 15=9+6, 7=5+2, and 3=3+0. It implies the
following numerical coincidence.

\proposition\label{magic} For any non-singular real cubic surface,
$$
h=\NN^+_\R,\ \ e=\NN^-_\R. \quad \quad\qed
$$
\endproposition

\subsection{A direct proof of Proposition \ref{magic} }
Consider a real line $\ell$ on a real non-singular projective
cubic surface $X$. Choose projective coordinates $x,y,u,v$ in
$P^3$ so that the line $\ell$ is given by equations $x=y=0$. Then
the defining polynomial of $X$ has the form
$f=u^2L_{11}+2uvL_{12}+v^2L_{22}+uQ_1+vQ_2+C$, where $L_{ij}$,
$Q_i$, and $C$ are homogeneous polynomials in $x$ and $y$ of
degree, respectively, one, two, and three. According to Segre's
definitions, the line $\ell$ is elliptic (respectively,
hyperbolic) if the quadratic form $L_{12}^2-L_{11}L_{22}$ is
definite (respectively, indefinite).

\lemma A real line contributes $1$ into $\NN_\R^+$ if the line is
hyperbolic, and to  $\NN_\R^-$ if elliptic.
\endlemma

\proof A line $\ell'$ in a neighborhood of $\ell\in G_2(\R^4)$
hits the coordinate projective planes $v=0$ and $u=0$ at some points
$[x_1\!:y_1\!:\!1\!:\!0]$ and $[x_2\!:y_2\!:\!0\!:\!1]$, and thus,
$\ell'$ can be given in a parametric form as $(u,v)\mapsto
u(x_1,y_1,1,0)+v(x_2,y_2,0,1)$.
 The value of section $s_f(\ell')$ is defined by the restriction of
$f$ to $\ell'$, i.e, by substitution $x=ux_1+vx_2$ and
$y=uy_1+vy_2$ in the polynomial $f$.
 Letting $L_{ij}=l_{ij}x+m_{ij}y$ for all $i,j$, we obtain
$s_f(\ell')=u^2[l_{11}(ux_1+vx_2)+m_{11}(uy_1+vy_2)]+\dots$. It is
straightforward to check, following the polar correspondence, that
our orientation of the total space of $Sym^3(\tau^*_{2,4})$ agrees
with the order $x_1,x_2,y_1,y_2$ of the local coordinates in
$G_2(\R^4)$ and the frame $u^3,u^2v,uv^2,v^3$ that gives a local
trivialization of $Sym^3(\tau^*_{2,4})$. The Jacobi matrix of
$s_f$ in these coordinates is as follows
$$
\left[\matrix
  l_{11}&0&m_{11}&0\\
  2l_{12}&l_{11}&2m_{12}&m_{11}\\
  l_{22}&2l_{12}&m_{22}&2m_{12}\\
  0&l_{22}&0&m_{22}\\
\endmatrix\right]
$$
Now, it remains to notice that the determinant of this matrix,
whose sign is the index $I_{\ell'}^e(s_f)$, is opposite to the
discriminant of the binary quadratic form $L_{12}^2-L_{11}L_{22}$.
\endproof

\section{Black Holes}
\subsection{Congruences}
As was disclosed by Gr\"un\-berg and Moree  \cite{GM}, the
residues of $\NN_\C$ modulo $2^q$ form a $2^q$-periodic sequence
for all $q\ge1$. For the modulo
$2^q$ residues of $\NN^e_\R=(2n-1)!!$
the $2^q$-periodicity is evident.  Up to $q=2$ the
two sequences of residues coincide, $\NN_\C=\NN^e_\R\mod 4$. The
relation between $8$-residues is less trivial; the period is
$1,1,3,3,5,5,7,7$ for $\NN_\C$ and $1,1,3,7,1,1,3,7$ for
$\NN^e_\R$. What can be general rules and conceptual explanations?

In the case of cubic surfaces another congruence holds,
$\NN_\R=\NN_\C \mod 4$. However, for any fixed dimension $n\ge3$,
the range of $\NN_\R$ includes all the odd integers in a certain
interval from some $\NN_\R^{\roman{min}}\ge\NN_\R^e$ to some
$\NN_\R^{\roman{max}}\le\NN_\C$.
 In a sense, this phenomenon
reflects the difference in the Galois groups of the corresponding
enumerative problems (starting from 3-folds the Galois group is
the full symmetric group, see \cite{H}).

\subsection{Welschinger invariants}\label{Welschinger}  There is a method, introduced by  J.Y.~Welschin\-ger \cite{We},
of counting real rational curves in a given homology class on
arbitrary real rational surfaces and certain real 3-folds. As was
observed by J.Solomon \cite{Sol}, in the case of quintic 3-folds
the number $\NN_\R^e=15$  can be seen as the Welschinger invariant
of the generator in the second homology group of the
complexification of the quintic. To get
 $\NN_\R^e=3$ in the case of cubic surfaces, one should take the sum of the Welschinger invariants over all homology classes of
projective degree $1$. Though, the invariants should be taken not
with the signs traditionally applied in the case of rational
surfaces, but with signs copied from that used in the case of
3-folds (just consider the $\operatorname{Pin}^-$-structure
$q:H_1(X_\R;\Z/2)\to \Z/4$ inherited by the real locus $X_\R$ of a
non-singular cubic surface from the ambient real projective space
and define the sign of a line $l$ to be equal to $i^{q(l)-1}$). As
an additional similarity,  let us  point that Welschinger
invariants tend to have asymptotic and arithmetic properties
analogous to that of the Euler numbers of real lines treated in
this note (see, for example, \cite{IKS}). Is there a general
theory of real enumerative invariants that bring all these results
together?

\subsection{ Hypersurfaces of other degrees}
If $X$ is a real hypersurface of odd degree $d<2n-1$ in
$P^{n+1}(\R)$, the dimension of the variety $F_\R(X)$ of real
lines contained in $X$ is  $\ge 2n-1-d$, and is  equal to $2n-1-d$
for generic $X$; in particular, it is always nonempty (see
\cite{DM}).  What is about abundance? Is there a lower bound for
the total Betti number, $b_*(F_\R(X))$ comparable in logarithmic
scale with  $b_*(F_\C(X))$?


\Refs\widestnumber\key{ABCD}

\ref{BVV}
 \by  W.~Barth, A.~Van de Ven
 \paper Fano varieties of lines on hypersurfaces
 \jour Arch. Math. (Basel)
 \vol 31
 \yr 1978/79
\pages 96--104
\endref\label{BVV}

\ref{DM}
\by  O.~Debarre, L.~Manivel
\paper Sur les intersections compl\`etes r\' eelles
\jour CRAS
 \vol 331 (S\'erie I)
\yr 2000
\pages 887--992
\endref\label{DM}

\ref{GM}
 \by D.B.~Gr\"unberg, P.~Moree
 \paper Sequences of enumerative geometry: congruences and asymptotics (with Appendix by Don Zagier)
 \jour Experiment. Math.
 \vol 17
 \yr 2008
\pages 409--426
\endref\label{GM}

\ref{H}
 \by J.~Harris
 \paper Galois groups of enumerative problems
 \jour Duke Math.~J.
 \vol 46
 \yr 1979
\pages 685--724
\endref\label{H}


\ref{IKS}
 \by I.~Itenberg, V.~Kharlamov, E.~Shustin
 \paper Welschinger invariants of real Del Pezzo surfaces of degree $\ge 3$
 \jour arXiv:1108.3369
\pages 32 pages
\endref\label{IKS}


\ref{OT}
\by    C.~Okonek,  A.~Teleman
\paper Intrinsic signs and lower bounds in real algebraic geometry
\jour arXiv:1112.3851
\pages 17 pages
\endref\label{OT}


\ref{Sc}
\by L.Schl\"afli
 \paper An attempt to determine the twenty-seven lines upon a surface of the third order, and to
divide such surfaces into species in reference to the reality of the lines upon the surface
 \jour  Quart. J. Pure Appl. Math.
 \vol 2
 \yr 1858
\pages 110--120
\endref\label{Sch}

\ref{Se}
 \by B.~Segre
\book The Non-Singular Cubic Surfaces. A new method of Investigation with Special Reference to
Questions of Reality
\bookinfo Oxford Univ. Press, London
 \yr 1942
\endref\label{Segre}

\ref{So} \by  J.P.~Solomon
\paper   Intersection theory on the
moduli space of holomorphic curves with Lagrangian boundary
conditions \jour arXiv:math/0606429 \pages 79 pages
\endref\label{Sol}


\ref{We}
\by J.-Y.~Welschinger
\paper Invariants of real rational symplectic 4-manifolds and lower
bounds in real enumerative geometry
\jour C. R. Acad. Sci. Paris, S\'er. I
\vol  336
\yr 2003
\pages 341--344\endref\label{We}


\endRefs
\enddocument